\documentclass[pdflatex,sn-mathphys-num]{sn-jnl}

\usepackage{graphicx}
\usepackage{amsmath,amssymb,amsfonts}
\usepackage[title]{appendix}
\usepackage{xcolor}
\usepackage{booktabs}
\usepackage{colortbl}
\newcommand{\sig}[1]{\cellcolor{red!12}#1}
\usepackage{algorithm}
\usepackage{algorithmicx}
\usepackage{algpseudocode}

\begin{document}

\title[New spherical codes by tuning an annealing solver]{New spherical codes by tuning an annealing solver}

\author*[1]{\fnm{Ilmari} \sur{Vahteristo}}\email{i.vahteristo@gmail.com}

\affil*[1]{\orgname{Independent Researcher}, \orgaddress{\city{Helsinki}, \country{Finland}}, ORCID: 0009-0005-3323-3286}

\abstract{
The spherical code problem asks how to place $N$ points on a unit sphere in $D$ dimensions to maximize the smallest distance between two points.
It is a non-smooth and non-convex max-min problem, where few problem sizes have known optimal solutions, some have constructions, but most still only have numerically obtained solutions.
We use a simple solver, with seven parameters, that anneals a LogSumExp surrogate, minimizes with L-BFGS, applies a random perturbation between anneals, and repeats with no acceptance test until a fixed limit is reached.
We search the solver parameters for every instance with $N \le 1024$, $3 \le D \le 32$, and $N > 2D$, a total of 29,670 instances, using 32 parameter trials. We then run the solver for 60 seconds with the best parameters for each instance, without utilizing any symmetries or existing solutions.
After enforcing monotonicity of the results and numerically polishing the solutions, we improve the best known code in 4,577 instances. When a single set of solver parameters is used for all instances, with one 60-second episode, the number of improvements is still 3,649. From our records, we find 24 codes with an exact algebraic cosine.
Among our ablations, we find that 1) the repeated anneal contributes most to improvements, 2) at our $32 \times 60$-second budget only fitting one parameter per-instance is better than fitting all seven, and 3) good parameters can be predicted from the instance size.
We provide the solver code, the fitted parameters, and the codes.
}



\maketitle

\section{Introduction}
\label{sec:intro}

The optimal placement of points on a 3D unit sphere to maximally separate the two closest points on the surface of the sphere is known as the Tammes problem~\cite{tammes1930}.
In higher dimensions, the problem is often called the spherical code problem~\cite{delsarte1977,conway1999}. The same mathematical problem has come up in various fields, such as biology~\cite{tammes1930}, chemistry and physics~\cite{atiyah2003}, coding theory and error correction~\cite{conway1999,ericson2001}.

Mathematically, the spherical code problem asks for $N$ unit vectors in $\mathbb{R}^D$ that maximize the minimum pairwise angle of the vectors. With the vectors written as a matrix $X$ of unit rows, this minimizes the largest pairwise cosine,
\begin{equation}
f(X)=\max_{i<j} x_i^\top x_j .
\label{eq:maxcos}
\end{equation}

Lai et al. note that the problem in three dimensions can act as a test system for global optimization techniques~\cite{lai2023}. At scale, and at variable dimensions, the problem has two more properties that make it a useful test. First, the instance size affects the geometry of the solution and thus the optimization landscape and the best parameters for a solver. Second, it allows us to study solver behavior and the problem geometry across different instance sizes and dimensions.

We approach the problem with one numerical solver that anneals a smooth LogSumExp surrogate of the objective, minimizes it with L-BFGS~\cite{byrd1995}, and perturbs the points with Student-$t$ noise. After every perturbation the solver returns to a smooth surrogate and anneals again, with no acceptance test. The solver is parameterized by seven parameters that govern its behavior.
We first fit one global parameter set with CMA-ES on a small sample of instances, and run the solver with those parameters on every instance. We then coarsely fit the solver parameters for every instance independently to maximize performance and analyze how the best parameters vary with $N$ and $D$. The solver is explained in Section~\ref{sec:episode}, the parameter fitting procedure in Section~\ref{sec:cma}, and the results in Section~\ref{sec:full-domain-runs}.

In Section~\ref{sec:res-params} we analyze the fitted solver parameters and Section~\ref{sec:freshseeds} analyzes how they depend on the number of points $N$ and the dimension $D$ and how predictable good parameters are.

We do ablation studies and comparisons. We ablate the repeated annealing and perturbation steps of our solver, study the effects of per-instance fitting, and compare different ways of using extended computational budgets with our solver (Section~\ref{sec:res-ablation}). We also compare against two recent numerical solvers in Section~\ref{sec:res-compare}.

We provide a reproducible episode environment, the solver implementation, our fitted parameters, and exact algebraic cosines for 24 codes. We also archive three tables of codes: the global table, the per-instance table, and a combined table that holds, for each instance, the best code found by any of our runs. These are found in the Data Availability section.

\subsection{The public table}
A well-maintained \emph{public table} records the best known code for every instance with $N \le 1024$, $3 \le D \le 32$, and $N > 2D$~\cite{cohn_spherical}. It gives an external reference value for every instance, and aggregates decades of work on finding optimal spherical codes. The entries (reported in cosine of the angle) in the public table are reported with an accuracy of $10^{-12}$, and any improvement over that numerical accuracy is accepted as an improvement to the table. Since the table is updated regularly, we compare our results against a snapshot of the table taken on 2026-08-02, unless otherwise specified.

The table excludes $N \le 2D$, since in 1955, Rankin proved that for $N \le D+1$ the regular simplex is optimal, and for $D+2 \le N \le 2D$ any $N$ vertices of a regular orthoplex attain the optimal angle $\pi/2$~\cite{rankin1955}.

In the rest of the domain, whether a given spherical code is optimal for an instance, characterized by the number of points $N$ and the dimension $D$, is proven only for 26 instances~\cite{cohn_spherical,schutte1951,danzer1986,hars1986,robinson1961,musin2012strong,musin2015tammes,boroczky1978,bachoc2009,delaat2024d4,levenshtein1979,odlyzko1979,cohn2024exact}.

Of the 29,670 instances, 16,265 have an explicit entry in our snapshot of the public table. The instances which do not have an explicit entry are not missing. An omitted entry at some $(N, D)$ instance means that a solution with a larger angle than one at $(N+k, D)$ (where $k$ is the smallest positive integer such that $(N+k, D)$ has an entry) has not been found. This is because a configuration at $(N+k, D)$ gives a configuration for $(N, D)$ that is at least as good, because the angle cannot decrease by removing a point. The table omits these to stay compact, and improvements to them are treated like improvements to any other instance. We call these omitted entries \emph{implied entries}.

\subsection{Related work}
Tables of known best configurations are maintained by Cohn~\cite{cohn_spherical}, and earlier ones were compiled by Sloane~\cite{sloane_tables}. The records in the public table come from algebraic constructions and from numerical search. Ericson and Zinoviev obtain codes by algebraic construction~\cite{ericson2001}, and Ganzhinov and \"Osterg\aa rd search for codes with symmetries inside shells of integer lattices~\cite{ganzhinov2024}.

The published numerical solvers mostly share one pattern. They replace the non-smooth objective \eqref{eq:maxcos} with a smooth surrogate, minimize it locally, and add a global strategy to escape local optima.
Nurmela approximates the packing problem in three to five dimensions by a repulsion energy, minimized by local and stochastic global optimization~\cite{nurmela1995}. He combines it with simulated annealing and tabu search, which minimize a penalty on the pairs closer than a target distance. Kottwitz minimizes an inverse-power repulsion energy in three dimensions at 15 exponents, doubling from 80 to 1,310,720~\cite{kottwitz1991}, and the procedure is repeated from random starts. Jay et al. minimize a Riesz interaction energy, optionally under symmetry constraints~\cite{jay2026}. In three dimensions, Lai et al. instead minimize a penalty on the pairs closer than a target distance~\cite{lai2023}. From a local optimum they perturb all points and optimize locally several times, and vary the number of trials and the perturbation size with the number of iterations since the last improvement. Hars does not use a surrogate in three dimensions, but instead minimizes a bound $d$ subject to $x_i^\top x_j \le d$ over the pairs with sequential quadratic programming~\cite{hars_tammes,hars_midsize}

LogSumExp is used as a smooth approximation of the maximum in finite minimax problems~\cite{pee2011}. Continuation methods for global optimization first minimize a smoothed function that has fewer local minimizers, then reduce the smoothing step by step, and start each minimization from the previous result~\cite{more1997}. The exponent doubling of Kottwitz and the anneal of our solver follow the same idea.

An iterated local search is a loop that perturbs a local optimum, optimizes locally again, and decides by an acceptance test whether to continue from the new optimum~\cite{lourenco2003}. Basin hopping is such a loop with a Metropolis acceptance test~\cite{wales1997,leary1999}. A variant that accepts only improvements~\cite{leary1999}, called monotonic basin hopping, has been applied to packing equal and unequal circles in a circular container~\cite{grosso2010} and to packing disks in a square~\cite{addis2008}. Our solver is an iterated local search in which the local search is one anneal. It always continues from the most recent configuration, which is known as the random walk acceptance criterion~\cite{lourenco2003}.

Sometimes the objective is reframed as a different type of problem. For example, PackingStar treats kissing configurations as a matrix-completion game over cosine entries, and trains a reinforcement-learning system to play it~\cite{ma2025kissing}.

Large language models (LLMs) have also seen much use in the design of spherical codes~\cite{novikov2025alphaevolve,bianchi2026harnessing}. ImprovEvolve uses an LLM to write the code for three steps, an initialization, a local improvement and a perturbation. A basin-hopping loop then runs them and keeps only improvements. The local step of ImprovEvolve is also a LogSumExp continuation, over the geodesic distances, optimized with L-BFGS~\cite{kravatskiy2026}.

Estimating solver parameters automatically has been studied as well. SMAC~\cite{hutter2011smac} and irace~\cite{lopez2016irace} search for one parameter set that performs well over a set of instances, and ISAC~\cite{kadioglu2010isac} one per cluster. We differ from these approaches in that we fit parameters for every instance individually to maximize the performance, improve records, and analyze how solver parameters should vary with $N$ and $D$.

\section{Methodology}
\label{sec:method}

This section describes the solver and the episode environment, how we fit the seven parameters, how we post-optimize the results, and how we recognize exact codes. The experiments and their results are in Section~\ref{sec:results}.

\subsection{The solver and its episode}
\label{sec:episode}
An \emph{episode} is a reproducible run of our solver on one instance, with a fixed time limit. The episode defines a controlled environment for the solver, where the starting configuration and the perturbation seed are fixed, which gives us reproducible results. To further ensure reproducibility, each episode is pinned to a single CPU core, so it shares only memory bandwidth and caches with other episodes running in parallel.

An episode essentially maps an instance, a parameter set, a perturbation seed, and a time limit to an angle.
More on the hardware, parallelism, and reproducibility of an episode is in Section~\ref{sec:hardware}.

Fitting (Section~\ref{sec:cma}) searches for the parameters that maximize the value of an episode (per-instance), or a normalized average of episodes across instances (global).

Our solver uses LogSumExp \eqref{eq:lse} as the surrogate objective, raises its hardness progressively, and applies a perturbation at the maximum hardness, after which it drops the hardness significantly, and restarts the anneal from the perturbed configuration. Our solver returns the best found configuration, but it does not have an acceptance test inside an episode; it always perturbs the current configuration, and repeats the anneal from it. Four of the seven parameters control the anneal, two the perturbation, and one the stopping tolerance of the optimization. Algorithm~\ref{alg:episode} shows how one episode progresses with our solver, and Figure~\ref{fig:traces} illustrates the behavior and the solution's response.

\subsubsection{Starting configuration and perturbation seed}
\label{sec:epinit}
To make comparing improvements across different instances fair, and not biased by poor starting configurations, we initialize the solver with a Sobol sample~\cite{sobol1967}, greedy farthest point sampling~\cite{gonzalez1985}, and one L-BFGS optimization call of Eq.~\eqref{eq:lse} at a fixed $\beta = 100$ and stopping tolerance $\epsilon \approx 2.2\times10^{-9}$, which is the default in the SciPy implementation of L-BFGS~\cite{virtanen2020scipy}.
The Sobol sample is a uniform sequence that covers the sphere more evenly than random sampling, and farthest-point selection ensures that the points are well separated even in the initialization. The starting configuration is shared across all episodes of the same instance, so the comparison between runs is fair. Since the initialization is precomputed and cached once, we do not count it toward the episode's time limit. Computing these starting positions for every instance took a total of 25 minutes on our hardware (Section~\ref{sec:hardware}).

The single L-BFGS call alone increases the median minimum angle by $11.8^\circ$ over the Sobol-farthest-point initialization, but still falls $2.60^\circ$ short of the public table and does not improve any instance of the public table.

The perturbation seed fixes the random draws of an episode. Two parameter sets with the same $\nu$ draw the same raw noise, and $\sigma$ only rescales the noise. Changing $\nu$ changes the perturbations. Runs with the same parameters and seed are reproducible up to hardware, floating point and wall-clock timing.

\begin{algorithm}[t]
\caption{One episode}
\label{alg:episode}
\begin{algorithmic}[1]
\Require Starting configuration $X_0$; time limit $T$; perturbation seed;
         solver parameters $\beta_{\min}, r, s, \rho, \sigma, \nu, \epsilon$
\State $X \gets X_0$; $X^* \gets X$;
       $\beta \gets \beta_{\min}$; $\beta_{\max} \gets r\,\beta_{\min}$
\While{wall clock $< T$}
  \State $X \gets \textsc{L-BFGS}(\mathcal{L}_\beta, X, \epsilon)$
         \Comment{optimize at fixed $\beta$}
  \If{$\min\text{-angle}(X) > \min\text{-angle}(X^*)$}
    \State $X^* \gets X$
           \Comment{new best; continue from $X$}
  \EndIf
  \If{$\beta < \beta_{\max}$}
    \State $\beta \gets \min(\beta \cdot r^{1/s},\, \beta_{\max})$
           \Comment{increase beta}
  \Else
    \State $X \gets \textsc{Perturb}(X, \sigma, \nu)$
           \Comment{perturb at $\beta_{\max}$}
    \State $\beta \gets \beta_{\max}^{1-\rho}\,\beta_{\min}^{\rho}$
           \Comment{roll back}
  \EndIf
\EndWhile
\State \Return $X^*$
\end{algorithmic}
\end{algorithm}

\subsubsection{Optimization}
\label{sec:optimization}
The objective \eqref{eq:maxcos} is non-smooth, so we anneal it with a smooth surrogate, the LogSumExp $\mathcal{L}_\beta$ \eqref{eq:lse}. The hardness parameter $\beta$ is the smoothing parameter of the surrogate.
The LogSumExp surrogate
\begin{equation}
\mathcal{L}_\beta(X)=\frac{1}{\beta}\log\sum_{i<j}\exp\bigl(\beta\, x_i^\top x_j\bigr)
\label{eq:lse}
\end{equation}
converges to \eqref{eq:maxcos} as $\beta \to \infty$.
For $M = N(N-1)/2$ pairs, $f(X) \le \mathcal{L}_\beta(X) \le f(X) + \log M/\beta$, so at $N = 1024$ and $\beta = 10^6$ the surrogate is within $1.3\times10^{-5}$ of the largest cosine. It is proven to be nearly optimal among smooth approximations of a maximum~\cite{samakhoana2025elementary}.

The only parameter that governs the optimization is the stopping tolerance $\epsilon$ (SciPy's \texttt{ftol}, $\frac{f_k - f_{k+1}}{\max{(|f_k|, |f_{k+1}|,1)}} < \epsilon$). The possible stopping criteria of the L-BFGS implementation we use are the stopping tolerance $\epsilon$, L-BFGS line search failure (flat loss landscape), the episode time limit, or a maximum of 20,000 iterations (which we never observe).

At each non-perturb step, we choose $\beta$ using the annealing parameters (Section~\ref{sec:annealing}) and minimize the surrogate with L-BFGS. We normalize each point when evaluating the loss and account for this normalization in the analytic gradient. The value and the gradient are computed in 64-bit floating point accuracy.

We do not store an $N \times N$ matrix of inner products directly, but rather compute and use it in blocks to save memory.
The points are processed in blocks of 128. For each block we compute the inner products of its points with all $N$ points, a $128 \times N$ array, and exponentiate it in place.
The array gives the block's part of the sum in \eqref{eq:lse} and gradient contributions for the block's points, and is then overwritten by the next block without reallocation. Computing and immediately using the exponential weights block by block, without storing the full matrix, uses a memory-saving approach similar to that of FlashAttention~\cite{dao2022flashattention}.

The exponentials of block $k$ are taken relative to the block's largest inner product $m_k$, as $\exp\bigl(\beta(x_i^\top x_j - m_k)\bigr)$, and after the last block, the sum and the gradients of each block are multiplied by $\exp\bigl(\beta(m_k - m)\bigr)$, where $m$ is the largest $m_k$.
This rescaling to combine blocks computed with different maxima is similar to online softmax~\cite{milakov2018online}.

Tested at $N = 1024$, $D = 32$, a simple kernel that allocates new $N \times N$ matrices on every call uses about 44 MiB, but only 1.8 MiB with the block implementation. On a single core ours completes 3 to 7 times as many L-BFGS iterations per second, compared to the simpler alternative, across problem sizes.
This matters for speed, and because concurrent episodes share memory bandwidth and cache (Section~\ref{sec:hardware}).

\begin{figure}[t]
\centering
\includegraphics[width=\textwidth]{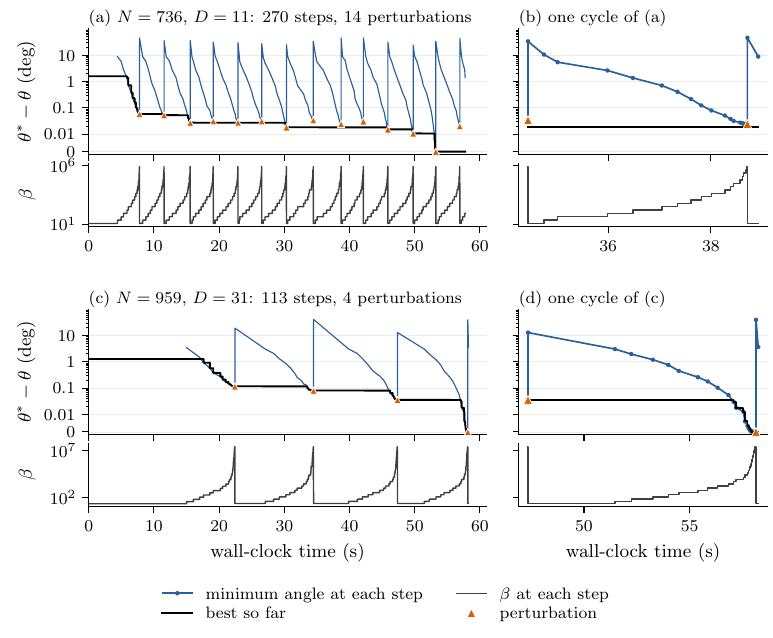}
\caption{Two examples of the solver behavior within 60-second episodes. The solver executes a sawtooth pattern with the $\beta$ parameter, and perturbs the points with fixed perturbation parameters every time $\beta_{\max}$ is reached and lowers $\beta$.
(a, c) The whole episode. (b, d) One cycle of it, from just before one perturbation to just after the next.
Upper panels: the shortfall of the minimum angle from the best angle of the episode, $\theta^* - \theta$, on a scale that is linear below $0.01^\circ$ and logarithmic above. Blue: at each step. Black: the best so far. Orange triangles: perturbations. Lower panels: $\beta$ at each step.}
\label{fig:traces}
\end{figure}

\subsubsection{Annealing}
\label{sec:annealing}

Annealing means progressively tightening the hardness of the objective. Our anneal is a geometric ladder, starting from $\beta_{\min}$ and multiplying by a constant factor $r^{1/s}$ after each optimization step, until $\beta = \beta_{\max} = r\,\beta_{\min}$. After reaching $\beta_{\max}$, the solver perturbs the points and lowers the hardness to $\beta = \beta_{\max}^{1-\rho}\,\beta_{\min}^\rho$. We call one optimization call at one $\beta$ a \emph{step}, and one climb of a ladder, ending at $\beta_{\max}$, an \emph{anneal}.

The four parameters $\beta_{\min}$, $r$, $s$, and $\rho$ govern the annealing behavior. $\beta_{\min}$ sets where the first anneal starts. A small $\beta_{\min}$ starts on a smooth objective where many angles matter at once, but which follows the true objective less closely. The span $r$ sets the largest $\beta$, and thus how closely the last calls approximate the maximum. The step count $s$ sets how many steps the annealing takes to get to $\beta_{\max}$: too few might freeze the configuration in poor local optima, too many take long to finish an anneal. The rollback $\rho$ sets how far down the annealing ladder the solver returns after each perturbation; only the first anneal is always guaranteed to start from $\beta_{\min}$.

Figure~\ref{fig:traces} shows the solver running on two
instances. The same 60-second time limit allows 270 steps at $N = 736$, $D = 11$ and 113 at $N = 959$, $D = 31$. We can also see the best found solution improving over time, as more anneals are executed, though irregularly, and with diminishing returns.

\subsubsection{Perturbation}
\label{sec:perturbation}
The objective is highly non-convex, so the solver is prone to getting stuck in local optima at high $\beta$ values. One way to move the solver to other local optima is to perturb the points with isotropic Student-$t$ noise: Each point receives an independent noise vector, projected onto its tangent plane and applied along the geodesic, so the points stay on the sphere's surface.
The scale $\sigma$, in radians, sets the size of the move, and $\nu$ sets how it is spread over the points. For large $\nu$, every point moves about $\sigma$. For small $\nu$, most points move a few $\sigma$, but a few move very far: at $\nu = 0.5$, the median point moves about $2\sigma$ but about one point in fifty moves more than $1{,}000\sigma$.

Mathematically, for each point $x_i$ independently, draw $z_i \sim \mathcal{N}(0, I_D)$
and $w_i \sim \chi^2_\nu$, and set
\begin{align}
v_i &= \sqrt{\nu / w_i}\;\bigl(I - x_i x_i^\top\bigr)\, z_i,
\label{eq:kicknoise}\\
x_i' &= \cos(\theta_i)\, x_i \;+\; \sin(\theta_i)\, v_i / \lVert v_i\rVert,
\label{eq:kickstep}
\end{align}
where the arc length is
$\theta_i = \sigma \lVert v_i\rVert / \sqrt{D-1}$. Every perturbation is followed by a rollback, so the perturbed configuration is optimized again from a smoother objective. Even a small perturbation can therefore lead to a different local optimum.

\subsection{Parameter search}
\label{sec:cma}
The parameter search runs in two stages. A global fit finds one parameter set on a small sample of instances. A per-instance fit refines that set separately on every instance. The refinement serves three purposes: more improved instances, a lower bound on what the solver can reach, and a study of how good parameters vary with $N$ and $D$.

To find the best parameters, globally and per-instance, we use CMA-ES~\cite{hansen2016cma}, a derivative-free evolutionary search algorithm, that is suitable for our seven-dimensional continuous parameter space where the obtained angle is a noisy and highly non-convex function of the parameters: A slightly different decision in the episode leaves the solver in a different local optimum, and that difference compounds, so two parameter sets that are close together can return very different angles.

CMA-ES works in generations. At each generation it proposes a small population, which we call \emph{candidates}, scores every candidate, and updates its search distribution from the ranking of those scores. The time allowed to optimize a set of parameters is determined by two factors: The episode length, and the number of episodes (which is the population size times the number of generations). A sweep of the two factors at a fixed compute budget is in Appendix~\ref{app:budget}.

A per-instance fit simply runs one episode per candidate and scores the candidate by the returned angle.
The global fit scores a candidate on a whole panel, a fixed subset of instances defined in Section~\ref{sec:results}, with one episode per instance.
Typical angles differ widely between instances: 500 points in three dimensions are much closer to each other than 200 points in twenty dimensions.
We therefore score each episode by how much it raises the angle over the starting configuration (Section~\ref{sec:epinit}), $\theta - \theta_0$, and rank candidates by the mean of $\log(\theta - \theta_0)$ over the panel.

Each parameter's search space is bounded.
This makes search more efficient, because the search does not waste time where the solver's behavior hardly changes between parameter changes. For example, increasing the degrees of freedom of the Student-$t$ distribution from 100 to 1,000 has little effect.

We search six of the parameters on a logarithmic scale and the rollback $\rho$ on a logit scale. The initial CMA-ES step size is a quarter of each parameter's range on that scale, so for most parameters the bounds lie two step sizes from the start. The stopping tolerance $\epsilon$ and the step count $s$ get a 1.75 times smaller step, because the angle is more sensitive to them.

The global fit starts at the midpoint of each parameter's range, while the per-instance fit starts from the global parameters. The parameters a fit returns are those of its best evaluated candidate, not the final CMA-ES mean.

Inside an episode two limits of our implementation apply. $\beta$ is capped at $10^9$. A parameter set with $r\,\beta_{\min} > 10^9$ therefore never reaches $\beta_{\max}$ and never perturbs. It simply runs one anneal from $\beta_{\min}$ and keeps optimizing at $\beta=10^9$ when it reaches it. This allows the solver to express a pure annealing behavior, without repeats or perturbations.
The perturbation scale $\sigma$ is kept between $2\times10^{-6}$ and $3$ times $N^{-1/(D-1)}$, which scales as the angle between neighboring points.

The ranges, midpoints, and initial step sizes for all parameters are summarized in Table~\ref{tab:searchspace}.

\begin{table}[t]
\centering\setlength{\tabcolsep}{4pt}
\caption{The seven solver parameters: the search range, the value the global fit starts from, and the values one initial step size below and above that start.}
\label{tab:searchspace}
\begin{tabular}{@{}lllll@{}}
\toprule
 & role & range & start & start $\pm$ one step size \\
\midrule
$\beta_{\min}$ & start of the anneal & 0.05 to 200 & 3.2 & 0.40 to 25 \\
$r$ & span $\beta_{\max}/\beta_{\min}$ & $10^2$ to $10^{10}$ & $10^6$ & $10^4$ to $10^8$ \\
$s$ & steps from $\beta_{\min}$ to $\beta_{\max}$ & 2 to $10^4$ & 141 & 43 to 470 \\
$\sigma$ & perturbation scale (rad) & $10^{-6}$ to 1.5 & $1.2\times10^{-3}$ & $3.5\times10^{-5}$ to 0.043 \\
$\nu$ & perturbation degrees of freedom & 0.5 to 100 & 7.1 & 1.9 to 27 \\
$\rho$ & rollback after a perturbation & 0.001 to 0.999 & 0.5 & 0.12 to 0.88 \\
$\epsilon$ & stopping tolerance & $10^{-13}$ to $10^{-6}$ & $3.2\times10^{-10}$ & $3.2\times10^{-11}$ to $3.2\times10^{-9}$ \\
\bottomrule
\end{tabular}
\end{table}

\subsection{Post-optimization}
\label{sec:postopt}

Our post-optimization pipeline is intended to enforce monotonicity, and to polish the best found solutions to the accuracy used in the public table. The post-optimization steps are only applied when a full table of results is available, and not on subsets of the full domain where we experiment, because the post-optimization pipeline requires successive $N$ values to enforce monotonicity and comparisons of results do not require polishing to a high numerical accuracy.

\textbf{1) We enforce monotonicity} in each of our tables. We go through every $D$ value, keep it fixed, and start at $N=1024$. If deleting a point from the solution for $N$ is better than our solution for $N-1$, we update the solution for $N-1$ and move to the next smaller $N$.

\textbf{2) We polish the solutions to a higher accuracy.} After enforcing monotonicity, we optimize each configuration by purely annealing $\beta$ in a fixed ladder starting from $10^2$, multiplying by 10 at each step, and stopping when the improvement between successive steps is $<10^{-14}$ in cosine, or $\beta>10^{13}$. The L-BFGS call at each $\beta$ value uses $\epsilon=10^{-18}$, and additionally we use a maximum of 2,000 iterations. This step aims to achieve the precision required by the public table, so that the reported angles cannot be overtaken by trivial polishing, but does not guarantee it.

\textbf{3) We alternate between enforcing monotonicity and polishing the solutions.}
After enforcing monotonicity, the solutions might not be tight anymore (removing a point leaves an empty space). Thus, we do a polishing step to bring the modified solutions closer to their true optima. Conversely, polishing the solutions might violate monotonicity as the positions of points are adjusted and tightened. Therefore, we alternate between the two for four rounds.

\subsection{Finding exact codes}
\label{sec:exactmethod}

A numerical configuration can sit close to a code whose largest cosine is an exact algebraic number. We call such a code an \emph{exact code}. We look for exact codes among the configurations that beat the public table, to ``snap'' our numerically optimized solutions to exact values when possible. A contact pair is a pair of points whose cosine equals the largest cosine $t$, and a rattler is a point that belongs to no contact pair. We keep the configurations that have more contact pairs than free coordinates. Then, on the contact pairs we solve $x_i^\top x_j = t$, with every $\lVert x_i \rVert = 1$, for the coordinates and $t$ together, to at least 200 digits. We then recognize $t$ as the root of an integer polynomial with PSLQ~\cite{ferguson1999pslq}, and express the entries of the Gram matrix in the field that root generates, or in a quadratic extension of it.
We certify the results by an $LDL^\top$ factorization of the Gram matrix in the number field, which shows that the points exist in $\mathbb{R}^D$. Then interval arithmetic on a rational interval that contains $t$ confirms that every inner product outside the contact pairs is below $t$.

The contact equations do not determine the positions of rattlers, and they can move slightly without changing $t$. In the certificate we fix each rattler in place and verify that all of its cosines are below $t$.
A certified code can sometimes take more points without increasing $t$, so after certifying a code we attempt to fill in such points. We keep all the points that have contacts fixed in place, remove the $k$ rattlers, and place $k+1$ points, by minimizing the largest cosine that involves any of the $k+1$ points from random starting positions. If this cosine is at most $t$, we keep the larger code and repeat with one more point. We then certify the code with the added points, as described above.

Finding the exact codes is the least reproducible part of our research. This is because it is not strictly in the scope of this paper, but exact codes are nonetheless very important for studying spherical codes. All the exact codes discussed in this paper have been initially found by our solver, and recognized as promising candidates using standard methods, however some instances required additional steps to finalize the exact codes.

\subsection{Hardware}
\label{sec:hardware}

All measurements ran on the Roihu supercomputer at CSC~\cite{csc_roihu_docs}. Each node has two 192-core AMD EPYC 9965 processors, 384 cores in all. Each episode runs single-threaded, pinned to one core, with the pinning balanced over 8 NUMA memory domains.

Contention in the node can still happen due to shared resources like DRAM access bandwidth. That contention reduces the number of steps an episode completes in a fixed time. We chose to use 192 concurrent episodes as our standard setting, which all runs that depend on timing use.

To determine the noise floor in our experiments, we measured how much the results change when a run is repeated on the same node architecture, but on different nodes. Across repetitions on a fixed sample of 1,024 instances, 92\% of instances reached the same cosine within $10^{-12}$, and the 95th percentile of absolute differences was $2.3 \cdot 10^{-5}$ in cosine.

\subsection{Use of generative AI}
\label{sec:genai}
I used Anthropic's Claude models through Claude Code and the web interface in this work. I used them to generate parts of the code from instructions, debug and test scripts, double-check claims and results, discuss and review the work, proofread the text, and research papers.
I have verified the claims of the paper and take full responsibility for its content.

\section{Results}
\label{sec:results}

In this section, we go through results from our large scale experiments, test alternative approaches, and perform ablations.
We first go through the results of our two large scale runs, where we first fit parameters, and then use the parameters on every instance of the domain. We focus on analyzing the results from our per-instance fitting, and the global parameters act mostly as a comparison. After analyzing these results, we decompose what contributes to our results.

The experiments and ablations use representative subsets of instances. We use two different panels: \textbf{r144} with 144 instances, and \textbf{r1024} with 1,024 instances, drawn randomly and independently.
The instances, starting configurations, and perturbation seeds for each panel are fixed and shared across all experiments. One perturbation seed per instance is used by every fit, and we call it the \emph{fitting seed}. Any other perturbation seed is a \emph{fresh seed}, which changes the random draws of the perturbations.

An instance is \emph{improved} when its cosine is smaller than in the table we compare against by more than $10^{-12}$, \emph{lost} when it is larger by more than $10^{-12}$, and \emph{tied} otherwise. An instance improved against the public table is a new \emph{record}.

We compare results by counting the instances improved, lost and tied, because the counts are agnostic to the scale differences between instances. The angles in our domain vary from $7^\circ$ to $89^\circ$, so whether an improvement of $0.01^\circ$ is considered a large gain or not depends on the instance.

\subsection{Full domain runs}
\label{sec:full-domain-runs}

The analysis here focuses on comparing the results of per-instance fitting against the public table, and the global table. Both the global and per-instance runs use 60-second episodes, due to our tests in Appendix~\ref{app:eplen-vs-trials}.

The global parameters are fitted on r144. We optimize the parameters by running 32 generations of CMA-ES, with a population size of 16. We then take the single set of parameters that gave the best score, and replay it on every instance with one 60-second episode.

The results from the replay are then put through the four rounds of post-optimization of Section~\ref{sec:postopt}. We call the resulting, post-optimized table the \emph{global table}.

The \emph{per-instance table} follows the same process, but instead of fitting one set of parameters on r144, we fit the parameters for each instance. Each instance uses 32 episodes to try parameters, divided into 4 generations of 8 candidates. From these trials, the parameters for that instance are the parameters that gave the largest angle. We do not use the best angle value obtained during the 32 trials, and instead we replay them, like the global approach does.

Table~\ref{tab:records} summarizes the global and per-instance tables compared to each other, before and after post-optimization. The global table sets 3,649 records, of which 32\% come from post-optimization, while the per-instance table sets 4,577 records, of which 19\% come from post-optimization. The post-optimization adds many records, especially to the global table, because many of the solutions have only partly converged, and increasing $\beta$ and decreasing $\epsilon$ makes easy gains.

The per-instance table reaches a record on 928 new instances compared to the global table, with 32 times the episodes. What record counts might hide is that against the global table, the per-instance table improves on 85\% of instances, 25,203 out of 29,670.

We also provide a third table of records in Appendix~\ref{app:combined-table}, which we generally do not discuss in this paper, since it is not the focus of this paper, and does not match the standards of reproducibility we otherwise exercise. This table contains, for each instance, the best solution obtained from the global and per-instance tables, as well as all our other experiments with our solver during this research, and improvements found during a search for the exact cosines we identify.

\begin{table}[t]
\centering
\caption{The full-domain runs against the public table.}
\label{tab:records}
\begin{tabular}{@{}lrr@{}}
\toprule
all 29,670 instances & per-instance & global \\
\midrule
\multicolumn{3}{@{}l}{\emph{before post-optimization}} \\
\quad instances improved & 3,708 & 2,455 \\
\quad mean gain on improved, $^\circ$ & 0.067 & 0.045 \\
\multicolumn{3}{@{}l}{\emph{after post-optimization}} \\
\quad instances improved & 4,577 & 3,649 \\
\quad mean gain on improved, $^\circ$ & 0.104 & 0.082 \\
\quad improved over the other table & 25,203 & 1,988 \\
\bottomrule
\end{tabular}
\end{table}

\begin{figure}[t]
\centering
\includegraphics[width=0.55\textwidth]{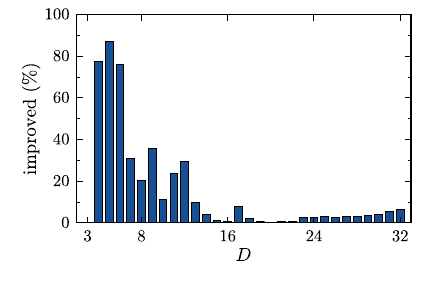}
\caption{Share of each dimension's instances that the per-instance table improves.}
\label{fig:perdim}
\end{figure}

\begin{figure}[t]
\centering
\includegraphics[width=\textwidth]{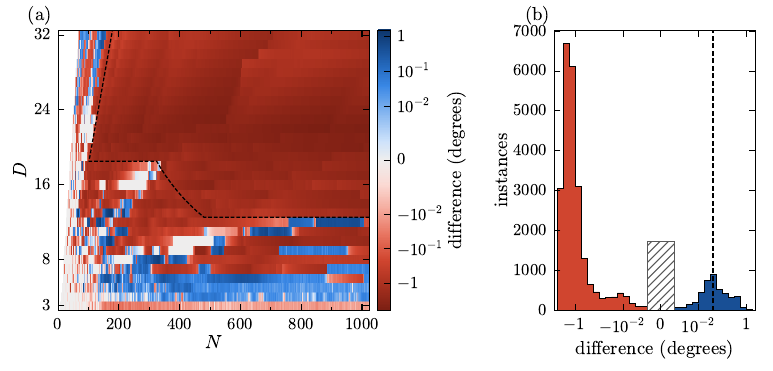}
\caption{The difference between the per-instance table and the best known code on all 29,670 instances, in degrees, on a scale that is logarithmic on both sides of zero. Blue: ours is larger. Red: ours is smaller. (a) The difference over $N$ and $D$. The dashed line marks the bounds beyond which we set no records, $ND = 6{,}000$ for $13 \le D \le 18$ and $N = 5.5D$ for $D \ge 19$. (b) The distribution of the difference. The hatched bar holds the 1,729 instances within $10^{-3}$ degrees, 921 of them tied. The dashed line marks the median improvement, $0.04^\circ$.}
\label{fig:records}
\end{figure}

\subsubsection{Per-instance table improvements}
\label{sec:per-instance-detail}
The per-instance table is a better target for analyzing the limitations of our solver, since it has better results and has tuned parameters for each instance.

Figure~\ref{fig:perdim} shows that the improvements are not equally spread across dimensions: We make many gains in small to mid dimensions from 4 to 13, few elsewhere, and none at $D=3$.

Neither the per-instance nor the global table sets or ties a record beyond two bounds in $N$ and $D$, as illustrated in Figure~\ref{fig:records}a.
For $13 \le D \le 18$, every record lies below $ND = 6{,}000$, and from $D \ge 19$, every record lies below $N = 5.5D$. 55\% of the domain is above these bounds, where the best known codes lie a median $2.0^\circ$ above ours.

The distribution of differences to the public table is shown in Figure~\ref{fig:records}b. The median improvement is $0.04^\circ$ and the largest $1.29^\circ$, with a long tail of improvements below $0.001^\circ$. 4,348 of the 4,577 improvements exceed $0.001^\circ$.

We take a new snapshot of the public table on 2026-09-23, with attribution data included, to compare against our per-instance table.
Comparing this snapshot with the previous version of the public table shows that 53\% of the instances were improved in those seven weeks. Most of the improvements are small, and 94\% of our per-instance records still hold in the new snapshot.

Table~\ref{tab:attrib} compares the improved instances of our per-instance table with the 2026-09-23 snapshot of the public table and its attribution data.
Of the 4,306 records still standing in the new snapshot, 164 improve entries whose attribution dates them to 2025 or earlier.

The oldest explicit entries that we improve are from Sloane's tables 1994--1996~\cite{sloane_tables}.
The oldest implied entries that the per-instance table improves had been derived from classical kissing configurations. The oldest implied entries we improve include $N = 180$ to $183$ in $D = 8$, which were previously subsets of the $\mathbf{E}_8$ root system found in 1873 by Korkine and Zolotareff~\cite{korkine1873}. The implied entries from $N = 250$ to $275$ in $D = 9$ were subsets from a code of 306 points found by Leech and Sloane in 1971~\cite{leech1971}.

\begin{table}[t]
\centering
\caption{Holder of the record our per-instance table improves over, computed from a 2026-09-23 snapshot of the public table. The table is also decomposed to explicit and implied entries. An implied entry is counted under the holder of the entry at the next larger $N$ that it follows from. Cohn is the table's maintainer.}
\label{tab:attrib}
\begin{tabular}{@{}lrrr@{}}
\toprule
record held by & explicit & implied & total \\
\midrule
Cohn, 2023--2026, flagged as needing more optimization & 3,099 & 38 & 3,137 \\
Jay, 2026~\cite{jay2026} & 653 & 48 & 701 \\
Borisov, 2026, unpublished & 203 & 10 & 213 \\
Cohn, 2023--2026, not flagged & 104 & 17 & 121 \\
Ericson and Zinoviev, 2001~\cite{ericson2001} & 1 & 68 & 69 \\
Leech and Sloane, 1971~\cite{leech1971} & 0 & 26 & 26 \\
Takhanov et al., 2026~\cite{takhanov2026} & 0 & 11 & 11 \\
anonymous, 2026 & 9 & 0 & 9 \\
Korkine and Zolotareff, 1873~\cite{korkine1873} & 0 & 6 & 6 \\
Li, 2026, unpublished & 1 & 5 & 6 \\
Ganzhinov and \"Osterg\aa rd, 2024~\cite{ganzhinov2024} & 1 & 2 & 3 \\
Sloane et al., 1994--1996~\cite{sloane_tables} & 2 & 0 & 2 \\
Ma et al., 2025~\cite{ma2025kissing} & 0 & 1 & 1 \\
Ericson and Zinoviev, 1995~\cite{ericson1995} & 0 & 1 & 1 \\
\midrule
total & 4,073 & 233 & 4,306 \\
\bottomrule
\end{tabular}
\end{table}

The table shows 73\% of our improvements are upon entries that had been flagged as ``needs more optimization'', which means that the configuration entry has ``not fully converged to a local optimum''~\cite{cohn_spherical}. Our search does not use the public table's configurations, so our records are not polishes over these configurations, unless our solver independently finds the same configuration and converges further.
The flag is on 70\% of the explicit entries in the public table, and our improvements over these are about 350$\times$ larger than the update between the snapshots: the median improvement on flagged entries was $1.3 \cdot 10^{-6}$ in cosine, while our median improvement on the same instances was $4.7 \cdot 10^{-4}$.

\subsubsection{Parameters}
\label{sec:res-params}

Table~\ref{tab:knobs} summarizes the fitted parameters of the 29,670 per-instance fits against the global parameters.

The fitted per-instance parameters are not the best parameters for each instance, since we only try 32 different parameters. However, some trends can be analyzed even if the parameters have not fully converged, but the data is noisy. A longer search on r144 is shown in Appendix~\ref{app:budget}.

The anneal rollback $\rho$ is near its maximum at 0.999 on almost every instance, and globally. Compared with the global fit, the per-instance parameters start the anneal lower, go higher, take more steps in an anneal, and stop each optimization call at a tighter tolerance. 

\begin{table}[!t]
\centering
\caption{The global parameters and a summary of the per-instance parameters. The per-instance columns summarize the 29,670 fits of the full-domain run.}
\label{tab:knobs}
\begin{tabular}{@{}lrrr@{}}
\toprule
 & global & per-instance median & 10th to 90th pct. \\
\midrule
$\beta_{\min}$ & 11.1 & 4.7 & 0.3 to 20 \\
$r$ & 7.5e4 & 5.5e5 & 1.4e4 to 8.6e7 \\
$s$ & 17 & 37 & 13 to 168 \\
$\sigma$ & 1.0e-3 & 1.0e-3 & 7.8e-6 to 0.092 \\
$\nu$ & 0.57 & 1.07 & 0.54 to 5.7 \\
$\rho$ & 0.999 & 0.997 & 0.972 to 0.999 \\
$\epsilon$ & 4.0e-8 & 7.5e-9 & 5.6e-10 to 5.2e-8 \\
\bottomrule
\end{tabular}
\end{table}

\begin{figure}[!t]
\centering
\includegraphics[width=\textwidth]{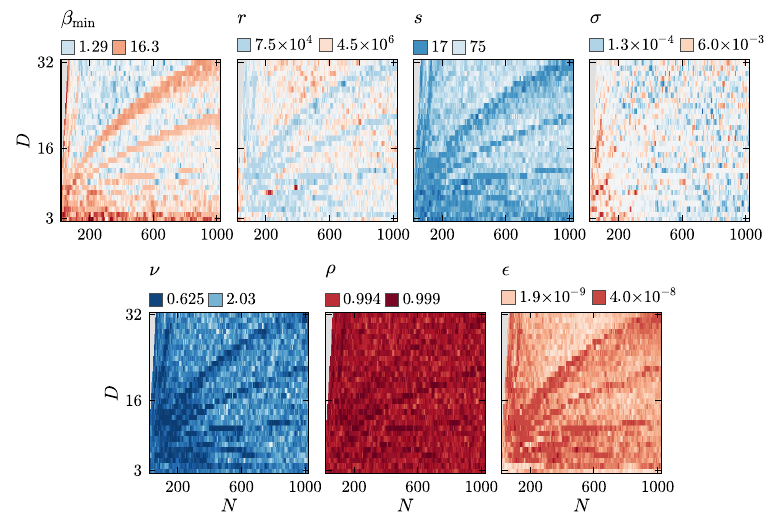}
\caption{The fitted value of each parameter over the domain, from the full-domain run, smoothed: each instance shows the median of its 16 nearest instances in $(\log N, \log D)$, which is the predictor used in Section~\ref{sec:freshseeds}. Red is above the search midpoint, white at it, and blue below. The two swatches above each panel give the color and the value at the 10th and 90th percentiles of that panel. The rollback is high nearly everywhere and the degrees of freedom low.}
\label{fig:fields}
\end{figure}

Figure~\ref{fig:fields} shows clear patterns in the preferred parameter values across the domain. The patterns show clear structure, which implies that the parameters should be predictable from instance size. We study this predictability in Section~\ref{sec:freshseeds}.

An interesting pattern is the streaks of slightly higher or lower parameter values that are shared across the parameters, and run across the plots indicating regions where the solver consistently prefers certain types of parameter settings. The streaks align with the geometric structure of the codes found by our solver, which is discussed in Appendix~\ref{app:antipodal}.

\subsection{Predicting good parameters}
\label{sec:freshseeds}

The per-instance parameters clearly outperform the global parameters when the perturbation seed is fixed, and the parameters are fitted against that seed. Figure~\ref{fig:fields} also shows clear structure in the parameters, which implies predictability. To decompose how much of the per-instance advantage is due to structural and predictable differences across instances rather than selection against a fitting seed, we test how much of the advantage can be obtained by predicting the parameters based on the instance's size.

We replay three parameter sets on r1024, at five fresh seeds: the global parameters, the per-instance parameters, and a \emph{predictor}. The predictor uses the per-instance fitted data and predicts the component-wise median of the fitted parameters of its 16 nearest other instances (excluding the instance itself) in $(\log N, \log D)$. It sees only $N$ and $D$, so it cannot be tuned to any instance's fitting seed. At fresh seeds each instance's angle is averaged over the five seeds.

If parameters did not matter, an instance would be improved or lost with equal probability, so we measure the advantage as improved minus lost instances.

Two measurements both indicate that nearly 75\% of the per-instance advantage over the global parameters is due to structural differences across $N$ and $D$, and not just favoritism to the fitting seed.

The per-instance parameters improve over the global parameters on 907 of the 1,024 instances and lose on 26 at the fitting seed (Table~\ref{tab:arms}). At fresh seeds they improve on 796 and lose on 153. The advantage that persists at fresh seeds is thus $(796 - 153) / (907 - 26) = 73\%$.

The predictor, on the other hand, improves over the global parameters on 833 instances and loses on 180 at the fitting seed. So the predictor has $(833 - 180) / (907 - 26) = 74\%$ of the per-instance advantage.

At fresh seeds the predictor improves over the global parameters on 885 instances and loses on 128. At fresh seeds the predictor beats the per-instance parameters on 558 instances and loses on 457, which is expected, since the per-instance parameters are partly tuned to the fitting seed.

\begin{table}[t]
\centering
\caption{Instances improved / lost / tied on r1024 between the global, predictor, and per-instance parameters. At fresh seeds each instance's angle is averaged over the five seeds.}
\label{tab:arms}
\begin{tabular}{@{}lr|r@{}}
\toprule
 & fitting seed & five fresh seeds \\
\midrule
per-instance vs global & \textbf{907} / 26 / 91 & \textbf{796} / 153 / 75 \\
predictor vs global & \textbf{833} / 180 / 11 & \textbf{885} / 128 / 11 \\
per-instance vs predictor & \textbf{697} / 318 / 9 & 457 / \textbf{558} / 9 \\
\bottomrule
\end{tabular}
\end{table}

\subsection{Comparison with two other solvers}
\label{sec:res-compare}

We ran the program SHS of Jay et al.~\cite{jay2026} and the best ImprovEvolve program~\cite{kravatskiy2026} on r144, with their authors' settings, with two exceptions that made SHS and ImprovEvolve better in our tests. We used random starting configurations for all solvers. We chose these two since SHS has multiple entries on the public table, ImprovEvolve also uses LogSumExp annealing, and because both are recent works.

All three programs ran on the same Roihu compute node, with the same number of workers. Both public solvers got 1,920 seconds per instance, the compute of a per-instance fit, while our solver only used 60 seconds with the global parameters. The details are in Appendix~\ref{app:compare}.

Table~\ref{tab:compare} shows the results, which indicate that at the public settings, and under our computational settings, our solver improves significantly over both solvers. Across eight seeds, our improved count ranges from 139 to 142 against SHS and from 98 to 108 against ImprovEvolve.

\begin{table}[t]
\centering
\caption{One 60-second episode with our global parameters against two other solvers on r144. SHS and ImprovEvolve both got 1,920 core-seconds per instance. The counts are comparisons against the median solution found by our solver at global parameters on each instance, over 8 fresh seeds. The median gain is the median over instances of our angle minus theirs.}
\label{tab:compare}
\begin{tabular}{@{}lrrrr@{}}
\toprule
 & \multicolumn{3}{c}{our episode is} & \\
\cmidrule(lr){2-4}
against & better & worse & equal & median gain \\
\midrule
SHS~\cite{jay2026} & 141 & 3 & 0 & $+0.093^\circ$ \\
ImprovEvolve~\cite{kravatskiy2026} & 102 & 42 & 0 & $+0.008^\circ$ \\
\bottomrule
\end{tabular}
\end{table}

\subsection{Ablations}
\label{sec:res-ablation}

In Table~\ref{tab:ablation} we compare the best angle each ablated method finds, with the per-instance table, before any post-optimization was applied. As a reference, rerunning the per-instance fit unchanged, except for the CMA-ES seed, improves 68 instances and loses 70 on r144, matching a coin flip, as it should.

\begin{table}[t]
\centering\setlength{\tabcolsep}{3pt}
\caption{Ablations of the per-instance fit, at the fitting seed. Each row uses the same amount of compute per instance as the per-instance fit: 32 episodes of 60 seconds. The improved/lost/tied counts compare the results with the per-instance table before post-optimization. $p$ is the two-sided sign test on the improved and lost counts, with ties excluded. Shaded: $p < 0.05$.}
\label{tab:ablation}
\begin{tabular}{@{}lllrr@{}}
\toprule
removed & replaced by & panel & improved / lost / tied & $p$ \\
\midrule
nothing & the same fit, rerun & r144 & 68 / 70 / 6 & 0.93 \\
\midrule
\multicolumn{5}{@{}l}{\emph{Is refinement with CMA-ES effective?}} \\
per-instance fit & one 1,920 s episode, global parameters & r144 & 46 / 97 / 1 & \sig{$2\times10^{-5}$} \\
per-instance fit & best of 32 seeds, global parameters & r144 & 45 / 99 / 0 & \sig{$8\times10^{-6}$} \\
CMA-ES & best of 32 random parameter sets & r144 & 49 / 90 / 5 & \sig{$6\times10^{-4}$} \\
\midrule
\multicolumn{5}{@{}l}{\emph{Which parameters need to be refined per instance?}} \\
refining $\beta_{\min}$ & $\beta_{\min}$ at its global value & r144 & 68 / 74 / 2 & 0.67 \\
refining $r$ & $r$ at its global value & r144 & 65 / 75 / 4 & 0.45 \\
refining $s$ & $s$ at its global value & r1024 & 477 / 506 / 41 & 0.37 \\
refining $\epsilon$ & $\epsilon$ at its global value & r144 & 41 / 95 / 8 & \sig{$4\times10^{-6}$} \\
refining $\sigma$ & $\sigma$ at its global value & r144 & 62 / 77 / 5 & 0.23 \\
refining $\nu$ & $\nu$ at its global value & r144 & 58 / 78 / 8 & 0.10 \\
refining $\rho$ & $\rho$ at its global value & r144 & 62 / 77 / 5 & 0.23 \\
refining all but $\epsilon$ & only $\epsilon$ fitted, the rest at global values & r1024 & 566 / 451 / 7 & \sig{$4\times10^{-4}$} \\
\midrule
\multicolumn{5}{@{}l}{\emph{Which parts of the solver are needed?}} \\
perturbations & no perturbations & r1024 & 454 / 569 / 1 & \sig{$4\times10^{-4}$} \\
repeated anneal & perturbations at $\beta_{\max}$ only ($\rho = 0$) & r144 & 24 / 119 / 1 & \sig{$2\times10^{-16}$} \\
no acceptance & perturbing the best configuration & r1024 & 497 / 523 / 4 & 0.43 \\
\bottomrule
\end{tabular}
\end{table}

\emph{Is refinement with CMA-ES effective?}
We spend the same 1,920 core-seconds per instance in three ways. The first runs the global parameters in one long 1,920-second episode. The second runs the global parameters with 32 different perturbation seeds and keeps the best. The third draws 32 parameter sets randomly from the initial CMA-ES distribution of the per-instance fit, and keeps the best. Each of the three loses about twice as many instances as it improves. Fitting the parameters per instance is therefore a better use of time than using the global parameters in a longer episode or more repetitions, and CMA-ES finds better parameters than random sampling.

\emph{Which parameters need to be refined per instance?} 
We hold one parameter at its global value and fit the other six per-instance as usual. Holding the tolerance $\epsilon$ fixed is most detrimental, and loses 95 instances and improves 41. Holding any of the other six fixed makes no detectable difference. Fitting $\epsilon$ alone, with the other six parameters at their global values, improves 566 instances and loses 451 on r1024. Thus, at our budget, fitting only $\epsilon$ improves more instances than refining all seven parameters. However, with eight times the compute, fitting only $\epsilon$ loses to the full seven-parameter fit of the same compute (Appendix~\ref{sec:res-deep}). So fitting $\epsilon$ is the most important step, but the other six parameters still clearly contribute to the overall performance. This is also evidenced in Figure~\ref{fig:fields}, which shows that neighboring instances have correlated parameter values in large streaks, which would not be the case if the other six parameters were not important.
Refining the step count $s$ and refining $\epsilon$ alone were tested on r1024 for higher statistical power.

\emph{Which parts of the solver are needed?}
Two ablations remove one part of the solver entirely, and fit the remaining parameters. Without perturbations, the solver still rolls back and anneals again, but never perturbs the configuration with noise.
With $\rho = 0$, the solver anneals once from $\beta_{\min}$ to $\beta_{\max}$, allowing a longer and more gradual annealing, but with no return to a smaller $\beta$ value. Removing the repeated anneal loses heavily, on 119 of 144 instances. Removing the perturbations has a smaller effect, but on r1024 it is a consistent and measurable effect.
The third ablation always perturbs the best configuration found so far, instead of perturbing the current configuration, which shows no difference.

\subsection{Exact codes}
\label{sec:res-exact}

We contribute exact algebraic cosines for 24 codes, which come from applying the procedure of Section~\ref{sec:exactmethod} to our records, not just those of the per-instance table. All the exact codes are listed in Table~\ref{tab:exact}, and released along with the code repository.

For most of the codes, every inner product lies in the field of $t$. At $N = 53$, $D = 24$; $N = 56$, $D = 25$; and $N = 64$, $D = 25$, the inner products lie in a quadratic extension of that field.

The $N=84$ code in $D=16$ only has a marginal improvement of $1.4\cdot10^{-12}$ in cosine over the public table. We take it to be the code of the public table, which our solver rediscovers, and we contribute its exact algebraic cosine.

Our attempt to enlarge the codes by filling in points was flawed. Henry Cohn found more points to place for three of them: one point of the $N = 68$ code in $D = 11$, three points of the $N = 332$ code in $D = 8$, and one point that enlarged a 151-point code in $D = 13$ to 152 points, which we later enlarged by two more points to $N = 154$. We report these enlarged codes.

\section{Discussion}
\label{sec:discussion}

Our solver sets multiple new records on the spherical code problem, over both recent and long-standing records. Notably, it does so without utilizing any existing solutions in the public table, or symmetries or existing constructions as is often done in the spherical code literature~\cite{nurmela1995,jay2026,ganzhinov2024,ericson2001,kravatskiy2026}. Some of the solutions we find with our solver have significant structure and symmetry, and their angle can be calculated exactly.

At our time limit the solver parameters seem to prefer exploring multiple local optima over converging to one, and our post-optimization step ensures that the coordinates converge more precisely.

The records our solver finds concentrate in low and middle dimensions, and the per-instance table has none in three dimensions, where most work has gone, and dedicated solvers exist~\cite{nurmela1995,kottwitz1991,lai2023,hars_tammes,hars_midsize}. Lai et al. alone use $50N$ (10 restarts with $5N$) minutes of computation per instance, corresponding to 167 hours at $N = 200$~\cite{lai2023}. In our combined table we find 39 improvements in $3D$ (Appendix~\ref{app:combined-table}).

Our solver sets no records beyond the bounds of Section~\ref{sec:per-instance-detail}, and this is a persistent limitation of our solver under our tests. The reason behind this is not clear, but in the region where we do not obtain any new records, over half of the entries (8,854 of 16,180) are implied entries that can be traced to a larger exact code.

Suitable parameters are not the same for all instances, and vary with the size of the spherical code. This can be seen by both the large improvements that per-instance parameters make over the global ones (Section~\ref{sec:full-domain-runs}), and by the amount of improvements that are obtainable by predicting the parameters from neighboring instances (Section~\ref{sec:freshseeds}). The instances have streaks of $N$ and $D$ where certain parameters are preferred. These streaks coincide with edges of regions where codes have antipodal structure (Appendix~\ref{app:antipodal}).

Among our ablations, we observe that removing the repeated anneal loses on the largest share of instances and that the perturbations are less critical, but still meaningful. In basin hopping~\cite{wales1997}, the perturbation is the only way to leave a local optimum and an acceptance test determines where the optimization continues from. Whether we repeat an anneal from the perturbed current solution, or from the perturbed best found solution does not make a difference. Repeating the anneal itself seems to allow the search to escape local optima. Our parameter search shows that it is preferable to repeat the anneal all the way from $\beta_{\min}$, rather than anywhere else between $\beta_{\min}$ and $\beta_{\max}$ (Section~\ref{sec:res-params}).

Parameter tuning ablations show that refining only the tolerance $\epsilon$ per instance, with the other six parameters at their global values, improves more instances than refining all seven on r1024 at our fitting budget. However, with 128 episodes of 120 seconds per instance, the full seven-parameter fit clearly wins against the fit of $\epsilon$ alone (Appendix~\ref{app:budget}).

Despite best efforts to control the repeatability of experiments, some variability remains since the episodes are time limited, and the performance of the computer can fluctuate between runs. The small number of episodes run for parameter fitting causes variability in the fitted parameters and the resulting performance. The search space of an ablated solver is contained within the full solver's search space, but in practice the parameters for an ablated solver might need to move a larger distance from the global parameters or span a smaller range to achieve comparable performance.

Further analysis of our spherical codes and their structure is out of the scope of this paper, and we hope future research will further explore our codes. We consider that especially the exact codes can be important, but also potentially ties with a known record, where our configuration might have a different structure than the known one, since the solution has been found independently.

Since a repeated annealing schedule is central to the solver's performance, future work could test reactive schedules that choose when to soften or tighten the objective. Since the tolerance $\epsilon$ is very important, it could similarly be made reactive.


\section{Conclusion}
\label{sec:conclusion}

We ran one simple solver for 60 seconds on every instance of the spherical code problem in $N \le 1024$, $3 \le D \le 32$, and $N > 2D$, without using symmetries or any existing codes.
With a single set of parameters it improves the public table on 3,649 instances after post-optimization steps.

With a 32$\times$ computational budget, we search for the best solver parameters using 32 trials, and after post-optimization, improve 85\% of the instances against a single parameter set, and find 4,577 records.

The solver parameters change with the instance size predictably, and about 75\% of the performance of 32 parameter fitting episodes can be obtained by predicting the parameters based on neighboring instances. The remaining 25\% seems to be parameter selection against a fixed perturbation seed.
At our budget of 32 episodes, fitting the stopping tolerance alone works better than fitting all seven parameters, but the other six parameters are useful with more compute.

We find that the solver works mainly because of the repeated anneal, and substituting the solver with a single, tighter annealing schedule loses a significant number of improvements.
When the solver returns to a smooth objective and anneals again, it can leave a local optimum even without perturbations. The perturbations are important, but minor compared to repeated annealing.

The solver sets records in one region of the domain only.
Inside two bounds in $N$ and $D$ it improves over a third of the public table, while beyond them, in 55\% of the domain, it sets no record at all, and eight times the compute does not close the gap.
We do not know why, but the codes of the public table in that region seem to have a structure that our solver does not find.

We leave further analysis of our own codes to future work. However, we publish all our records found during this research, and among our records we report 24 configurations that have an exact cosine. Since our search never uses the codes of the public table, our codes are found independently of them, and there may be new and interesting configurations among them.
All our codes, the solver, the fitted parameters, and the certificates of the exact codes are public, and we hope they are useful for that work.

\backmatter

\bmhead{Acknowledgments}
I thank Henry Cohn for maintaining the table of spherical codes, which provided the benchmark data and relevant literature, and for helpful discussion and feedback.
I wish to acknowledge CSC -- IT Center for Science, Finland, for generous computational resources on the Roihu supercomputer.

\section*{Declarations}

\bmhead{Funding}
This work did not receive any funding.

\bmhead{Competing interests}
There are no competing interests to declare.

\bmhead{Data availability}
\label{sec:data-availability}
The solver, the fitted parameters, and certificates for the 24 exact codes are available at \url{https://github.com/ilmari99/Spherical-code-solver} and archived at \url{https://doi.org/10.5281/zenodo.23270597}. The repository contains our solver code and its episodes reproduce ours on Roihu, up to wall-clock timing. 

We archive three tables of coordinates at \url{https://doi.org/10.5281/zenodo.23270586}. Each table holds one code for each of the 29,670 instances.

The global table and per-instance table are as explained in Section~\ref{sec:full-domain-runs}. The \emph{combined table} holds, for each instance, the best code found by any of our runs. We include it for completeness, and since it contains some new records, which should not be wasted. Its sources are the global and per-instance tables, the exact cosine search, the ablations, episode-length experiments, and some other experiments not documented in this paper. All the codes in the combined table are obtained from our solver, but might have used different time limits, parameters, post-optimization or steps done during exact cosine value search.

All results and comparisons pertaining to full-domain results presented in this paper use the global or per-instance tables only, with the only exception being the exact cosine results and Appendix~\ref{app:combined-table}. The combined table improves 5,191 instances of the public table (as of 2026-08-02). The codes we submitted to the public table at \url{https://spherical-codes.org} are those of the combined table.

\begin{appendices}
\makeatletter\@addtoreset{table}{section}\@addtoreset{figure}{section}\makeatother
\renewcommand{\theHtable}{\thesection.\arabic{table}}
\renewcommand{\theHfigure}{\thesection.\arabic{figure}}

\section{Compute budget}
\label{app:budget}

\subsection{Episode length vs more trials}
\label{app:eplen-vs-trials}
Episode length and the number of episodes for CMA-ES tuning both consume the same compute budget.
We run six per-instance fits on r144 at a constant 1,920 core-seconds per instance, with episode lengths of 3.75, 7.5, 15, 30, 60 and 120 seconds; 512 episodes at 3.75 seconds down to 16 at 120.
The search starts from global parameters, fitted for 15-second episodes. These are a separate fit from the 60-second global parameters of Section~\ref{sec:full-domain-runs}, which is why the 60-second row does not reproduce the paper's per-instance values.
The final parameters of each fit are replayed at their own episode length and Table~\ref{tab:budget} shows the comparisons with the per-instance table, before post-optimization was applied.

\begin{table}[!ht]
\centering
\caption{Longer episodes tend to win up to 60-second episodes, and the 120-second episodes are less effective at the same compute. Improved / lost compare each fit's parameters after 1,920 core-seconds per instance, with the per-instance table of the paper, before post-optimization. Steps and perturbations are medians over instances.}
\label{tab:budget}
\begin{tabular}{@{}rrrrr@{}}
\toprule
Time limit & episodes & improved / lost & steps & perturb. \\
\midrule
3.75 s & 512 & 22 / 122 & 28 & 0 \\
7.5 s & 256 & 31 / 113 & 50 & 1 \\
15 s & 128 & 37 / 107 & 94 & 1 \\
30 s & 64 & 47 / 97 & 210 & 2 \\
60 s & 32 & 55 / 89 & 540 & 8 \\
120 s & 16 & 39 / 105 & 1,022 & 30 \\
\bottomrule
\end{tabular}
\end{table}

The table shows the improved count rising with episode length up to 60 seconds, after which it declines since at 120 seconds per-instance, the CMA-ES only has 16 episodes, 2 generations to improve the parameters. This justifies our choice to use 60-second episodes for our experiments and runs, which had the same 1,920-second budget per-instance.

\subsection{More compute per instance}
\label{sec:res-deep}
When we increase the compute budget per-instance from 1,920 core-seconds to 15,360 core-seconds, corresponding to 120-second episodes and 128 episodes per instance, we improve 136 out of the 144 instances compared to the per-instance fits with 1,920 core-seconds.

Part of this gain could come from the longer episode alone. When we use the per-instance parameters that are tuned to 60-second episodes, but increase the episode length to 120 seconds, the new parameters tuned to 120-second episodes still improve on 106 instances. This indicates that most of the advantage is from tuning the parameters for 120 seconds, and from having more episodes to tune the parameters.

Table~\ref{tab:deep-params} shows the fitted parameters for the per-instance fits of the full-domain run and the fits with eight times the compute.
With longer episodes, the solver prefers to start the anneal softer and reach a higher $\beta$ ($\beta_{\max}$ $9.6\times10^{6}$ against $2.5\times10^{6}$). They take about twice as many steps per anneal and stop each optimization call at a three times tighter tolerance.

We also repeated the $\epsilon$-only fit of Section~\ref{sec:res-ablation} with 120-second episodes and 128 episodes per instance on r144, with the other six parameters at their global values. Compared with the above full seven-parameter fit at the same compute, fitting $\epsilon$ alone loses clearly, and loses on 106 of the 144 instances.

\begin{table}[t]
\centering
\caption{Fitted parameters on r144 for the per-instance fits of the full-domain run (60-second episodes, 32 episodes) and the fits with eight times the compute (120-second episodes, 128 episodes): the median and the 10th to 90th percentiles over the 144 instances.}
\label{tab:deep-params}
\begin{tabular}{@{}lr@{\hskip 2em}rr@{\hskip 2em}rr@{}}
\toprule
 & & \multicolumn{2}{c}{per-instance} & \multicolumn{2}{c}{eight times the compute} \\
\cmidrule(r{2em}){3-4}\cmidrule{5-6}
 & global & median & 10th to 90th pct. & median & 10th to 90th pct. \\
\midrule
$\beta_{\min}$ & 11.1 & 4.3 & 0.22 to 21 & 1.8 & 0.10 to 20 \\
$r$ & 7.5e4 & 7.7e5 & 2.5e4 to 1.2e8 & 7.6e6 & 7.5e4 to 7.2e8 \\
$s$ & 17 & 35 & 13 to 140 & 65 & 16 to 356 \\
$\sigma$ & 1.0e-3 & 8.8e-4 & 3.3e-6 to 0.052 & 1.0e-3 & 9.2e-6 to 0.12 \\
$\nu$ & 0.57 & 1.19 & 0.57 to 4.8 & 1.63 & 0.54 to 29 \\
$\rho$ & 0.999 & 0.998 & 0.973 to 0.999 & 0.996 & 0.946 to 0.999 \\
$\epsilon$ & 4.0e-8 & 6.4e-9 & 5.4e-10 to 4.0e-8 & 1.9e-9 & 6.3e-11 to 3.0e-8 \\
\bottomrule
\end{tabular}
\end{table}

91 of the 144 instances of r144 lie beyond the bounds identified in Section~\ref{sec:per-instance-detail}, where we set no records. There, the fits with eight times the compute improve 85 of the 91 over the per-instance table, yet do not significantly close the gap to the public table: the mean angle is still 2.99\% below the table's, against 3.01\% for the per-instance table.

\section{Exact algebraic cosines}
\label{app:exact}
The exact algebraic cosines are listed in Table~\ref{tab:exact} and certificates are released alongside the solver code.

\begin{sidewaystable}
\centering
\caption{Exact algebraic cosines certified.
Each cosine is a root of the stated irreducible minimal polynomial
(verified in exact rational arithmetic), and the exact Gram matrix has rank $D$ and is positive semidefinite.
The last column counts the points without contacts.}
\label{tab:exact}
\begin{tabular}{@{}rrllr@{}}
\toprule
$N$ & $D$ & cosine & minimal polynomial & points without contacts \\
\midrule
52 & 24 & 0.04545454545455 & {\scriptsize$22x - 1$} & 0 \\
64 & 25 & 0.05882352941176 & {\scriptsize$17x - 1$} & 0 \\
78 & 21 & 0.08965517241379 & {\scriptsize$145x - 13$} & 12 \\
166 & 9 & 0.44444444444444 & {\scriptsize$9x - 4$} & 22 \\
56 & 25 & 0.04619127393258 & {\scriptsize$79x^2 + 18x - 1$} & 0 \\
60 & 28 & 0.03970337884400 & {\scriptsize$584x^2 + 2x - 1$} & 0 \\
68 & 11 & 0.23181909492386 & {\scriptsize$79x^2 - 14x - 1$} & 4 \\
80 & 18 & 0.12141716201614 & {\scriptsize$31x^2 - 12x + 1$} & 0 \\
100 & 12 & 0.25984357433173 & {\scriptsize$167x^2 - 28x - 4$} & 15 \\
138 & 12 & 0.29911947447944 & {\scriptsize$49x^2 - 18x + 1$} & 10 \\
176 & 13 & 0.29911947447944 & {\scriptsize$49x^2 - 18x + 1$} & 16 \\
69 & 15 & 0.15528799292197 & {\scriptsize$193x^3 - 65x^2 - x + 1$} & 1 \\
84 & 16 & 0.15947325285064 & {\scriptsize$17856x^3 + 2064x^2 - 476x - 49$} & 0 \\
87 & 21 & 0.10320429547298 & {\scriptsize$168x^3 - 156x^2 + 24x - 1$} & 0 \\
60 & 26 & 0.04628367725636 & {\scriptsize$4416x^4 + 1568x^3 - 376x^2 - 8x + 1$} & 0 \\
78 & 12 & 0.23324577331977 & {\scriptsize$455568x^4 - 93824x^3 - 18708x^2 + 2448x + 289$} & 6 \\
154 & 13 & 0.29644926916925 & {\scriptsize$183325x^4 - 90460x^3 + 8458x^2 + 636x + 9$} & 46 \\
66 & 13 & 0.18007654564827 & {\scriptsize$197541x^6 - 1327980x^5 + 842815x^4 + 93800x^3 - 58725x^2 + 500x + 625$} & 2 \\
174 & 10 & 0.41148237435905 & {\scriptsize$9216x^6 - 9216x^5 - 9216x^4 + 7296x^3 - 992x^2 - 32x + 1$} & 39 \\
46 & 16 & 0.09371477347368 & {\scriptsize\begin{tabular}[t]{@{}l@{}}$793984x^8 + 2769888x^7 - 316860x^6 - 393992x^5$\\ $\quad{}- 62620x^4 - 1308x^3 + 460x^2 + 40x + 1$\end{tabular}} & 0 \\
53 & 24 & 0.04651491804364 & {\scriptsize\begin{tabular}[t]{@{}l@{}}$160403328x^8 + 74571840x^7 - 10873152x^6 - 3270912x^5$\\ $\quad{}+ 865608x^4 - 83220x^3 + 4048x^2 - 100x + 1$\end{tabular}} & 0 \\
70 & 20 & 0.09512792431150 & {\scriptsize\begin{tabular}[t]{@{}l@{}}$15691627956x^8 - 2820129980x^7 - 185860412x^6 + 44167092x^5$\\ $\quad{}+ 518453x^4 - 212164x^3 - 906x^2 + 364x + 1$\end{tabular}} & 0 \\
332 & 8 & 0.58620103383402 & {\scriptsize\begin{tabular}[t]{@{}l@{}}$49408x^8 - 2130944x^7 + 2977024x^6 + 7215872x^5$\\ $\quad{}- 4396192x^4 - 4177856x^3 + 976208x^2 + 727376x + 28561$\end{tabular}} & 2 \\
154 & 7 & 0.54601454613448 & {\scriptsize\begin{tabular}[t]{@{}l@{}}$125x^9 + 2745x^8 + 492x^7 - 2092x^6 - 1090x^5$\\ $\quad{}+ 94x^4 + 284x^3 + 84x^2 - 3x + 1$\end{tabular}} & 1 \\
\bottomrule
\end{tabular}
\end{sidewaystable}

\section{Comparison to two public solvers}
\label{app:compare}

Both programs ran on the same type of node as our episodes (Section~\ref{sec:hardware}), as subprocesses of 192 pinned workers. For each program we ran two configurations, and we report the one that wins more instances against our solver. We did not start the solvers from the public table, because the percentage of improvements over the public table is fairly low, so it would take more instances to go through to reach a count of improvements that is statistically significant to differentiate.

For \textbf{SHS}~\cite{jay2026}, we use commit 76a5083 of the authors' repository. The authors' own campaign scripts run the program from random points with a new seed for every run, and we do the same. We restart until the time runs out and keep the best finished run. Each run uses the stage settings of the authors' production runs, a Riesz energy followed by a dot-product stage. Their production runs end with an SQP stage, but without the SQP stage SHS finishes more runs and finds better codes. SHS finished at least one run on every instance, and a median of 8.

For \textbf{ImprovEvolve}~\cite{kravatskiy2026}, we use the best program from the authors' repository (\texttt{champion\_gen11\_199c6110.py}, gigaevo-core commit 2fccecc) under their own restart protocol. The protocol improves one starting code and then runs chains. Each chain perturbs the code at ten intensities from $1$ to $10^{-6}$ and runs the program's improvement step after each perturbation. A chain restarts from the best code, or from a new random configuration every fifth chain. The authors start from a code of the public table, but we start from one random configuration from ImprovEvolve's own random configuration generator. A variant that first improves ten random configurations was better on 44 instances and worse on 88, so we report the protocol with one start, which performed better. A step took a median of 39 seconds at $3{,}640 \le N \cdot D < 10{,}000$ and 117 seconds above that, and the program made a median of 30 steps per instance. In $8 \le D \le 16$, the dimensions ImprovEvolve was tested on in its paper, our median episode from random starting configurations improves 29 instances and loses 8.

Our global parameters were fitted on r144 and this risks giving us an advantage on r144. However, there were only 7 global parameters so there is not much capability to overfit. And regardless, the global parameters do not perform better on r144 than elsewhere: When replayed on all instances, they beat the public table on 6.3\% of r144 against an average of 8.3\% elsewhere. In addition, we handicap our solver by using the median of 8 fresh seeds instead of the fitting seeds, we use random starting configurations which the parameters are not tuned for, and a 32$\times$ smaller computational budget.

\section{The streaks and the structure of our codes}
\label{app:antipodal}

\begin{figure}[t]
\centering
\includegraphics[width=\textwidth]{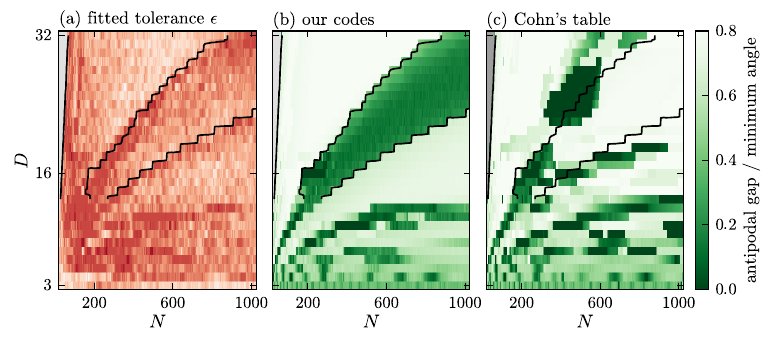}
\caption{The streaks and the structure of the codes. (a) The fitted tolerance $\epsilon$, as in Figure~\ref{fig:fields}. (b) The antipodal gap divided by the minimum angle for our codes, dark where near-antipodal. (c) The same for an April 2026 snapshot of the public table's codes, filled with implied entries where necessary. The black line, on every panel, outlines the range where our codes are near-antipodal: the ratio in (b), median-filtered over nine consecutive $N$, below 0.45, for $D \ge 13$.}
\label{fig:antipodal}
\end{figure}

The fitted parameters change along a few streaks in the domain (Figure~\ref{fig:antipodal}a). The streaks follow a change in the shape of codes our solver finds.

A code is antipodal when $-x$ is a point of the code for every point $x$. We measure how close our codes come to this with the antipodal gap: the median, over the points $x_i$, of the angle between $-x_i$ and the nearest point of the code. We consider a code \emph{near-antipodal} if this gap is below 0.45 times the minimum angle of the code.
Figure~\ref{fig:antipodal} illustrates the streaks and the antipodality of our codes and the public table's codes across $N$ and $D$.

The antipodal streaks in our codes match the changes in the fitted parameters and do not come from the per-instance fitting. The global table, made with one global parameter set, changes antipodal structure at the same places, within about 3\% in $N$.

The larger diagonal streak does not appear in the public table's codes, though there are some similarities, especially with the lower streak. In most of this range the codes of the public table are not near-antipodal by our measure, so the near-antipodal structure is a property of the codes our solver finds there.

\section{Combined table}
\label{app:combined-table}
This section gives the results of the combined table (see the Data availability statement) and compares it with the per-instance table.

The combined table improves more instances of the public table than the per-instance table does. It improves 5,191 instances of the 2026-08-02 snapshot, against 4,577 for the per-instance table. Of these, 4,979 still hold in the 2026-09-23 snapshot.

Unlike the per-instance table, the combined table has improvements in three dimensions: 39 instances at $N = 250$ and between $N = 449$ and $569$. The improvements in three dimensions were found by either further polishing our per-instance table, or by testing longer episodes, up to 120 seconds.

Table~\ref{tab:attrib-combined} lists the holders of the records that the combined table improves in the 2026-09-23 snapshot, as Table~\ref{tab:attrib} does for the per-instance table. Entries flagged as needing more optimization make up 69\% of them (3,452 of 4,979), against 73\% for the per-instance table.

\begin{table}[!hb]
\centering
\caption{Holder of the record our combined table improves over, computed from the 2026-09-23 snapshot of the public table. Explicit and implied entries are counted as in Table~\ref{tab:attrib}.}
\label{tab:attrib-combined}
\begin{tabular}{@{}lrrr@{}}
\toprule
record held by & explicit & implied & total \\
\midrule
Cohn, 2023--2026, flagged as needing more optimization & 3,404 & 48 & 3,452 \\
Jay, 2026~\cite{jay2026} & 836 & 57 & 893 \\
Cohn, 2023--2026, not flagged & 212 & 22 & 234 \\
Borisov, 2026, unpublished & 218 & 10 & 228 \\
Ericson and Zinoviev, 2001~\cite{ericson2001} & 2 & 81 & 83 \\
Leech and Sloane, 1971~\cite{leech1971} & 0 & 30 & 30 \\
Korkine and Zolotareff, 1873~\cite{korkine1873} & 0 & 14 & 14 \\
anonymous, 2026 & 10 & 1 & 11 \\
Takhanov et al., 2026~\cite{takhanov2026} & 0 & 11 & 11 \\
Li, 2026, unpublished & 2 & 8 & 10 \\
Ganzhinov and \"Osterg\aa rd, 2024~\cite{ganzhinov2024} & 1 & 4 & 5 \\
Sloane et al., 1994--1996~\cite{sloane_tables} & 3 & 0 & 3 \\
Ma et al., 2025~\cite{ma2025kissing} & 0 & 3 & 3 \\
Ericson and Zinoviev, 1995~\cite{ericson1995} & 0 & 1 & 1 \\
orthogonal union of two codes & 1 & 0 & 1 \\
\midrule
total & 4,689 & 290 & 4,979 \\
\bottomrule
\end{tabular}
\end{table}

\end{appendices}

\clearpage
\bibliography{samba_refs}

\end{document}